\documentclass[12pt,a4paper]{amsart}
\usepackage[latin1]{inputenc}
\usepackage{latexsym}
\usepackage{color,graphicx,shortvrb}
\usepackage{amsmath, amssymb}
\usepackage{amsfonts}
\usepackage[colorlinks, bookmarks=true]{hyperref}

\newtheorem{theorem}{Theorem}[section]
\newtheorem{lemma}[theorem]{Lemma}

\newtheorem{corollary}[theorem]{Corollary}

\theoremstyle{definition}

\author{J. M. Almira}
\title{A note on classical and p-adic Fr\'{e}chet functional equations with restrictions}

\begin{document}

\begin{abstract} Given $X,Y$ two $\mathbb{Q}$-vector spaces, and $f:X\to Y$, we study under which conditions on the sets $B_k\subseteq X$, $k=1,\cdots,s$, if  $\Delta_{h_1h_2\cdots h_s}f(x)=0$ for all $x\in X$ and $h_k\in B_k$, $k=1,2,\cdots,s$, then $\Delta_{h_1h_2\cdots h_s}f(x)=0$ for all $(x,h_1,\cdots,h_s)\in X^{s+1}$.
\end{abstract}

\maketitle

\markboth{J. M. Almira}{Fr\'{e}chet functional equation with restrictions}
\section{Introduction}
Let $X, Y$ be two $\mathbb{Q}$-vector spaces and let $f:X\to Y$. We say that $f$ satisfies Fr\'{e}chet's functional equation of order $s-1$ if \begin{equation}\label{fre}
\Delta_{h_1h_2\cdots h_s}f(x)=0 \ \ (x,h_1,h_2,\dots,h_s\in X),
\end{equation}
where $\Delta_hf(x)=f(x+h)-f(x)$ and $\Delta_{h_1h_2\cdots h_s}f(x)=\Delta_{h_1}\left(\Delta_{h_2\cdots h_s}f\right)(x)$, $s=2,3,\cdots$. In particular, after Fr\'{e}chet's seminal paper \cite{frechet}, the solutions of this equation are named ``polynomials'' by the Functional Equations community, since it is known that, under very mild regularity conditions on $f$, if $f:\mathbb{R}\to\mathbb{R}$ satisfies \eqref{fre}, then $f(x)=a_0+a_1x+\cdots a_{s-1}x^{s-1}$ for all $x\in\mathbb{R}$ and certain constants $a_i\in\mathbb{R}$. For example, in order to have this property, it is enough for $f$ being locally bounded \cite{frechet}, \cite{almira_antonio}, but there are stronger results \cite{ger1}, \cite{kuczma1}, \cite{mckiernan}. The general solutions of \eqref{fre} are characterized as functions of the form $f(x)=A_0+A_1(x)+\cdots+A_n(x)$, where $A_0$ is a constant and $A_k(x)=A^k(x,x,\cdots,x)$ for a certain $k$-additive symmetric function $A^k:X^k\to Y$ (we say that $A_k$ is a diagonalization of $A^k$). Furthermore, it is known that $f:X\to Y$ satisfies \eqref{fre} if and only if it satisfies the (apparently less restrictive) equation
\begin{equation}\label{frepasofijo}
\Delta_{h}^sf(x):=\sum_{k=0}^s\binom{s}{k}(-1)^{s-k}f(x+kh)=0 \ \ (x,h\in X).
\end{equation}
A proof of this fact follows directly from \cite[Theorem 9.2, p. 66]{czerwik}, where it is proved that the operators $\Delta_{h_1 h_2\cdots h_s}$ satisfy the equation
\begin{equation}\label{igualdad}
\Delta_{h_1\cdots h_s}f(x)=
\sum_{\epsilon_1,\dots,\epsilon_s=0}^1(-1)^{\epsilon_1+\cdots+\epsilon_s}
\Delta_{\alpha_{(\epsilon_1,\dots,\epsilon_s)}(h_1,\cdots,h_s)}^sf(x+\beta_{(\epsilon_1,\dots,\epsilon_s)}(h_1,\cdots,h_s)),
\end{equation}
where $\alpha_{(\epsilon_1,\dots,\epsilon_s)}(h_1,\cdots,h_s)=(-1)\sum_{r=1}^s\frac{\epsilon_rh_r}{r}$ and $\beta_{(\epsilon_1,\dots,\epsilon_s)}(h_1,\cdots,h_s)=\sum_{r=1}^s\epsilon_rh_r$.  (Note that $\Delta_{h}^sf(x)$ results from $\Delta_{h_1h_2\cdots h_s}f(x)$ by imposing $h_1=\cdots=h_s=h$).

Given $D\subseteq X$, the function $f:D\to Y$ is named a ``polynomial on $D$'' if $f$ satisfies \eqref{frepasofijo} for a certain $s\in\mathbb{N}$ and all $x,h\in X$ such that $\{x,x+h,\cdots,x+sh\}\subseteq D$. A natural problem, that has been solved by R. Ger \cite{ger}, is to study the conditions under which a polynomial $f$ on $D$ can be extended (and how to make this) to a polynomial on $X$. In this paper we study the analogous problem when we try to extend not the function $f$, which is assumed to be defined on all the space $X$, but the domain of validity of the equation \eqref{fre} for the steps $h_1,h_2,\cdots, h_s$. We prove our results for the special cases $X=\mathbb{R}$ and $X=\mathbb{Q}_p$, where $\mathbb{Q}_p$ denotes the field of $p$-adic numbers (which is, of course, a very special $\mathbb{Q}$-vector space). We also consider the case of $X$ being a topological vector space in a more general setting.

\section{The case $X=\mathbb{R}$}
In this section we assume that $f:\mathbb{R}\to Y$ for a certain $\mathbb{Q}$-vector space $Y$.

\begin{lemma}\label{lema} Let us assume that $I=(a,b)$ is a nonempty open interval of the real line. If $\Delta_hf(x)=0$ for all $x\in\mathbb{R}$ and $h\in I$ then $\Delta_hf(x)=0$ for all $x,h\in\mathbb{R}$.
\end{lemma}

\noindent \textbf{Proof.} We assume, with no loss of generality, that $0\leq a<b$. Obviously, $\Delta_hf(x)=0$ for all $x,h\in\mathbb{R}$ if and only if $f$ is a constant function. As a first step, we prove that $f$ is constant on an interval of the form $(\alpha,\infty)$ for a certain real number $\alpha$.

Given $x_0\in\mathbb{R}$, we know that $\Delta_hf(x_0)=f(x_0+h)-f(x_0)=0$ for all $h\in (a,b)$. Hence $f_{|x_0+I}=f(x_0)$. What is more, if $k\in\mathbb{N}$ and $h\in (a,b)$, then
\begin{eqnarray*}
0 &=& \Delta_hf(x_0+(k-1)h) = f(x_0+kh)-f(x_0+(k-1)h),\\
0 &=& \Delta_hf(x_0+(k-2)h) = f(x_0+(k-1)h)-f(x_0+(k-2)h),\\
\vdots & & \\
0 &=& \Delta_hf(x_0) = f(x_0+h)-f(x_0),
\end{eqnarray*}
so that $f(x_0)=f(x_0+kh)$. This implies that $f(x)$ is constant on the set $M=x_0+\bigcup_{k=1}^\infty(ka,kb)$. Now, the intervals $(ka,kb)$ and $((k+1)a,(k+1)b)$ overlap for all $k\geq k_0$ for a certain $k_0\in\mathbb{N}$, since $0\leq a/b<1$ and $k/(k+1)$ converges monotonically to $1$. Hence, there exists a real number $\alpha$ such that $f_{|(\alpha,\infty)}=f(x_0)$.

Let us prove that $f$ is constant on all the real line. We know that $f(x)=C$ for a certain constant $C$ and all $x\in (\alpha,\infty)$. On the other hand,  the same argument we used for the first part of the proof shows that, for all $x_1\leq \alpha$, if  $h\in (a,b)$ and $k\in\mathbb{N}$ is big enough, then $x_1+kh>\alpha$ and $f(x_1)=f(x_1+kh)=C$. {\hfill $\Box$}

\begin{theorem}\label{dos} Let us assume that $I_k=(a_k,b_k)$ are nonempty open intervals of the real line, $k=1,2,\cdots,s$.  If $\Delta_{h_1h_2\cdots h_s}f(x)=0$ for all $x\in\mathbb{R}$ and $h_k\in I_k$, $k=1,2,\cdots,s$, then $\Delta_{h_1h_2\cdots h_s}f(x)=0$ for all $(x,h_1,\cdots,h_s)\in\mathbb{R}^{s+1}$.
\end{theorem}

\noindent \textbf{Proof.}  By hypothesis, given $(h_2,h_3,\cdots,h_s)\in I_2\times I_3\times \cdots \times I_s$, the function $F(x)=\Delta_{h_2h_3\cdots h_s}f(x)$ satisfies the hypothesis of Lemma \ref{lema} with $a=a_1,b=b_1$. It follows that $\Delta_{h_1h_2\cdots h_s}f(x)=\Delta_{h_1}(\Delta_{h_2 \cdots h_s}f)(x)=0$ for all $x\in\mathbb{R}$ and $(h_1,h_2,\cdots,h_s)\in \mathbb{R}\times I_2\times\cdots\times I_s$. On the other hand, it is well known \cite[Corollary 9.1, p. 66]{czerwik} that, for any permutation $\sigma$ of the indices $\{1,2,\cdots,s\}$, we have that  $\Delta_{h_1h_2\cdots h_s}=\Delta_{h_{\sigma(1)} h_{\sigma(2)} \cdots h_{\sigma(s)}}$. Thus, if we take $\sigma_2=(12)$ (i.e., the transposition of the indices $\{1,2\}$, see \cite[p. 49]{jacobson} for the definition) and we apply the argument above to $\Delta_{h_{\sigma_2(1)} h_{\sigma_2(2)} \cdots h_{\sigma_2(s)}}$, we will get
$\Delta_{h_1h_2\cdots h_s}f(x)=0$ for all $x\in\mathbb{R}$ and $(h_1,h_2\cdots,h_s)\in \mathbb{R}\times \mathbb{R}\times I_3\times\cdots\times I_s$. The proof follows by repetition of the argument, taking into account the transpositions $\sigma_k=(1k)$, $k=3,4,\cdots,s$. {\hfill $\Box$}

\begin{theorem}\label{tres} Let $I=(-\delta,0)$ for a certain $\delta>0$. If $\Delta_{h}^sf(x)=0$ for all $x\in\mathbb{R}$ and $h\in I$, then $\Delta_{h_1h_2\cdots h_s}f(x)=0$ for all $(x,h_1,\cdots,h_s)\in\mathbb{R}^{s+1}$. An analogous result holds for $I=(0,\delta)$.
\end{theorem}
\noindent \textbf{Proof. } Let us assume that $-\delta/s\leq h_k\leq 0$ for $k=1,2,\cdots,s$. Take $(\epsilon_1,\cdots,\epsilon_s)\in\{0,1\}^s$ and set $\alpha_{(\epsilon_1,\dots,\epsilon_s)}(h_1,\cdots,h_s)=(-1)\sum_{r=1}^s\frac{\epsilon_rh_r}{r}$. Then
\[
0\leq \alpha_{(\epsilon_1,\epsilon_2,\dots,\epsilon_s)}(h_1,h_2,\cdots,h_s)=(-1)\sum_{r=1}^s\frac{\epsilon_rh_r}{r} \leq  \frac{1}{s}\left(\sum_{r=1}^s\epsilon_r\right)\delta\leq \delta.
\]
and, taking into account the equation \eqref{igualdad} above, it follows that we can use Theorem \ref{dos} with $I_k=(-\delta/s,0)$ for $k=1,\cdots,s$.

The last claim of the theorem follows from the relation that exists between the operators $\Delta_{-h}^s$ and $\Delta_h^s$:
\[
\Delta_{-h}^sf(x)=(-1)^s\Delta_{h}^sf(x-sh).
\]
{\hfill $\Box$}

\noindent \textbf{Remark 1.} Take $f(x)=x$ for $x\in\mathbb{Q}$ and $f(x)=x^2$ for $x\in\mathbb{R}\setminus \mathbb{Q}$. Then $\Delta_{h_1h_2h_3}f(x)=0$ for all $x\in\mathbb{R}$ and all $(h_1,h_2,h_3)\in\mathbb{Q}^3$. On the other hand, $f$ can not be a solution of Fr\'{e}chet's equation $\Delta_h^sf(x)=0$ for all $x\in\mathbb{R}$ and all $h\in\mathbb{R}$ for any $s\in\mathbb{N}$, since $f$ is not a polynomial function (in the ordinary sense) and $f$ is locally bounded. This shows that, in order to extend the set of validity of the parameters $h_i$ (as in Theorems \ref{dos},\ref{tres} above), it is quite natural to assume that the equation holds true for $h_i$ in a certain open set.
\medskip

\noindent \textbf{Remark 2.} Let $f:X\to Y$, where $X, Y$ are $\mathbb{Q}$-vector spaces and $X$ admits a topology $\tau_X$ with the property that all neighborhoods of the origin are ``naturally absorbent'' sets (we say that $B$ is naturally absorbent if for each $x\in X$ there exists $k\in\mathbb{N}$ such that $x\in kB=\{kz:z\in B\}$). If we assume that $\Delta_hf(x)=0$ for all
$x\in X$ and $h\in B$ for a certain neighborhood $B$ of the origin, then the arguments of the first part of the proof of Lemma \ref{lema} lead to the conclusion that $\Delta_hf(x)=0$ for all $x,h\in X$, since $X=\bigcup_{k\geq 1}kB$. It follows that Theorems \ref{dos}, \ref{tres} admit the following natural generalizations:

\begin{theorem}\label{cuatro} Let $f:X\to Y$, where $X, Y$ are $\mathbb{Q}$-vector spaces and $X$ admits a topology $\tau_X$ with the property that all neighborhoods of the origin are naturally absorbent sets. Let  $B$ be  a neighborhood of the origin.  If $\Delta_{h_1h_2\cdots h_s}f(x)=0$ for all $x\in X$ and $h_k\in B$, $k=1,2,\cdots,s$, then $\Delta_{h_1h_2\cdots h_s}f(x)=0$ for all $(x,h_1,\cdots,h_s)\in X^{s+1}$.
\end{theorem}

\begin{theorem}\label{cinco} Let $f:X\to Y$, where $X, Y$ are $\mathbb{Q}$-vector spaces and $X$ is a real normed vector space. Let  $B=B(0,\varepsilon)=\{x:\|x\|_X<\varepsilon\}$ be  an open ball centered in the origin of $X$. If $\Delta_{h}^sf(x)=0$ for all $x\in X$ and $h\in B$, then $\Delta_{h_1h_2\cdots h_s}f(x)=0$ for all $(x,h_1,\cdots,h_s)\in X^{s+1}$.
\end{theorem}

\section{The $p-$adic case} Let $\mathbb{Q}_p$ denote the field of $p$-adic numbers, which are expressions of the form
\begin{equation}\label{numero}
x=a_{m}p^{m}+a_{m+1}p^{m+1}+\cdots+a_0+a_1p+a_2p^2+ \cdots + a_np^n +\cdots =\sum_{n\geq m}a_np^n,
\end{equation}
where $m\in\mathbb{Z}$, $1\leq a_{m}\leq p-1$ and  $0\leq a_k\leq p-1$ for all $k>m$ (see, for example, \cite{gouvea}, \cite{robert} for the definition and basic properties of these numbers and their field extensions $\mathbb{C}_p$, $\Omega_p$).    Given $x$ as in \eqref{numero},  its $p$-adic absolute value is given by $|x|=p^{-m}$. The set of $p$-adic numbers with $|x|_p\leq 1$ is denoted by $\mathbb{Z}_p$. Obviously, $|x|_p\leq p^n$ if and only if $x\in p^{-n}\mathbb{Z}_p:=\{p^{-n}h: h\in\mathbb{Z}_p\}$. An important property of the absolute value $|\cdot|_p$ we will use is the following one:
\begin{equation}\label{propiedadnumeros}
 (|x|_p>|y|_p) \Rightarrow (|x+y|_p=|x|_p)
 \end{equation}
In this section we assume that $f:\mathbb{Q}_p\to Y$, where $Y$ is a $\mathbb{Q}$-vector space. Previous to any work about the extension of Fr\'{e}chet functional equation in this context, it would be appropriate to say something about the equation in the context of $p$-adic analysis.  In particular, we show that Fr\'{e}chet's original result has a natural extension to this new setting:

\begin{theorem}[p-adic version of Fr\'{e}chet's theorem]   \label{teo1}
Let $(\mathbb{K},|\cdot|_{\mathbb{K}})$ be a valued field such that $\mathbb{Q}_p\subseteq \mathbb{K}$ and the inclusion $\mathbb{Q}_p \hookrightarrow\mathbb{K}$ is continuous. Let us assume that $f:\mathbb{Q}_p\to \mathbb{K}$ is continuous. Then $f$ satisfies $\Delta_h^{n+1}f(x)=0$ for all $x,h\in\mathbb{Q}_p$ if and only if  $f(x)=a_0+\cdots+a_nx^n$ for certain constants $a_k\in \mathbb{K}$.
\end{theorem}

\noindent \textbf{Proof. } Assume $\Delta_h^{n+1}f(x)=0$ for all $x,h\in\mathbb{Q}_p$. Let $x_0,h_0\in\mathbb{Q}_p$ and let $p_0(t)\in\mathbb{K}[t]$ be the polynomial of degree $\leq n$ such that $f(x_0+kh_0)=p_0(x_0+kh_0)$ for all $k\in\{0,1,\cdots,n\}$ (this polynomial exists and it is unique, thanks to Lagrange's interpolation formula). Then
\begin{eqnarray*}
0 &=&
\Delta_{h_0}^{n+1}f(x_0)=\sum_{k=0}^n\binom{n+1}{k}(-1)^{n+1-k}f(x_0+kh_0)+f(x_0+(n+1)h_0)\\
&=& \sum_{k=0}^n\binom{n+1}{k}(-1)^{n+1-k}p_0(x_0+kh_0)+f(x_0+(n+1)h_0)\\
&=& -p_0(x_0+(n+1)h_0)+f(x_0+(n+1)h_0),
\end{eqnarray*}
since  $0=\Delta_{h_0}^{n+1}p(x_0)=\sum_{k=0}^{n+1}\binom{n+1}{k}(-1)^{n+1-k}p_0(x_0+kh_0)$. This means that $f(x_0+(n+1)h_0)=p_0(x_0+(n+1)h_0)$. In particular, $p_0=q$, where $q$ denotes the polynomial of degree $\leq n$ which interpolates $f$ at the nodes $\{x_0+kh_0\}_{k=1}^{n+1}$. This argument can be repeated (forward and backward) to prove that $p_0$ interpolates $f$ at all the nodes $x_0+h_0\mathbb{Z}$.

Let $m\in\mathbb{Z}$ and let us use the same kind of argument, taking $h_0^*=h_0/p^m$ instead of $h_0$. Then we get a polynomial $p_0^*$ of degree $\leq n$ such that $p_0^*$  interpolates $f$ at the nodes $x_0+\frac{h_0}{p^m}\mathbb{Z}$. Now, $p_0=p_0^*$ since the set
\[
\frac{h_0}{p^m}\mathbb{Z}\cap h_0\mathbb{Z}=\left\{
\begin{array}{cccccc}
h_0\mathbb{Z} &  & \text{if} &   m\geq 0 \\
\frac{h_0}{p^m}\mathbb{Z}  &  & \text{if} & m< 0\\
\end{array}
\right.
\]
is infinite. Thus, we have proved that $p_0$ interpolates $f$ at all the points of  $$\Gamma_{x_0,h_0}:=x_0+\bigcup_{m=-\infty}^\infty\frac{h_0}{p^m}\mathbb{Z}.$$
Now, $\Gamma_{0,1}$ is a dense subset of $\mathbb{Q}_p$, so that the continuity of $f$ implies that $f=p_0$ everywhere. {\hfill $\Box$}

\begin{corollary} [Local p-adic version of Fr\'{e}chet's theorem]
Let $(\mathbb{K},|\cdot|_{\mathbb{K}})$ be a valued field such that $\mathbb{Q}_p\subseteq \mathbb{K}$ and the inclusion $\mathbb{Q}_p \hookrightarrow\mathbb{K}$ is continuous. Let us assume that $f:\mathbb{Q}_p\to \mathbb{K}$ is continuous and let $N\in\mathbb{N}$. If  $f$ satisfies $\Delta_{p^N}^{n+1}f(x)=0$ for all $x\in\mathbb{Q}_p$, then for each $a\in\mathbb{Q}_p$ there exists constants $a_k\in \mathbb{K}$ such that $f(x)=a_0+\cdots+a_nx^n$ for all $x\in a+p^N\mathbb{Z}_p$.
\end{corollary}
\noindent \textbf{Proof. } Just repeat the argument of the first part of the proof of Theorem \ref{teo1}   with $x_0=a$ and $h_0=p^N$. This will show that there exists a polynomial of degree $\leq n$, $p_0\in \mathbb{K}[t]$ such that $(p_0)_{|a+p^N\mathbb{Z}}= f_{|a+p^N\mathbb{Z}}$. Now, $f$ is continuous and $a+p^N\mathbb{Z}$ is a dense subset of $a+p^N\mathbb{Z}_p$. {\hfill $\Box$}

It is well known that there are non constant locally constant functions $f:\mathbb{Q}_p\to Y$. In particular, the characteristic function associated to $\mathbb{Z}_p$, given by $\phi(x)=1$ for $x\in \mathbb{Z}_p$ and $\phi(x)=0$ otherwise, is continuous, non-constant, and $\Delta_h\phi(x)=0$ for all $x\in\mathbb{Q}_p$ and $h$ such that $|h|_p\leq 1$. Furthermore, $\Delta_{\frac{1}{p}}\phi(0)=\phi(p^{-1})-\phi(0)=-1\neq 0$. This is in contrast with the results we got for $X=\mathbb{R}$.

\begin{lemma}\label{lemapadic} Let $f:\mathbb{Q}_p\to Y$ and let $N\in\mathbb{Z}$, $a\in\mathbb{Q}_p$ be such that $\Delta_hf(x)=0$ for all $x\in\mathbb{Q}_p$ and all $h\in \mathbb{Q}_p\setminus (a+ p^{-N}\mathbb{Z}_p)$ then $\Delta_hf(x)=0$ for all $x,h\in\mathbb{Q}_p$.
\end{lemma}

\noindent \textbf{Proof.} We divide the proof in two cases:

 \noindent \textbf{Case 1: $a=0$.}
Given $x_0\in\mathbb{Q}_p$, we know that $\Delta_hf(x_0)=f(x_0+h)-f(x_0)=0$ for all $h\in \mathbb{Q}_p\setminus p^{-N}\mathbb{Z}_p$. Let us take $h\in p^{-N}\mathbb{Z}_p$. Then there exists $m\in \mathbb{N}$ such that  $\frac{h}{p^m}\not\in p^{-N}\mathbb{Z}_p$ and
\begin{eqnarray*}
0 &= & \Delta_{\frac{h}{p^m}}f(x_0+\frac{(p^m-1)}{p^m}h) = f(x_0+h)-f(x_0+\frac{(p^m-1)}{p^m}h),\\
0 &= & \Delta_{\frac{h}{p^m}}f(x_0+\frac{(p^m-2)}{p^m}h) = f(x_0+\frac{(p^m-1)}{p^m}h)-f(x_0+\frac{(p^m-2)}{p^m}h)=0,\\
\vdots & & \\
0 &= & \Delta_{\frac{h}{p^m}}f(x_0) = f(x_0+\frac{1}{p^m}h)-f(x_0),
\end{eqnarray*}
so that $f(x_0)=f(x_0+h)$.

\noindent \textbf{Case 2: $a\neq 0$.}   Let $k_0\in\mathbb{Z}$ be such that  $|a|_p=p^{k_0}$ and let $h\in \mathbb{Q}_p$ be such that $|h|_p=p^m$. If $m\neq k_0$ then $|h-a|_p=\max\{p^{k_0},p^m\}$, so that the imposition of $m\geq M=\max\{N,k_0\}+1$ implies that $|h-a|_p=p^m >p^N$. Hence $\mathbb{Q}_p\setminus p^{-M}\mathbb{Z}_p \subseteq \mathbb{Q}_p\setminus (a+ p^{-N}\mathbb{Z}_p)$ and we are again in Case 1. {\hfill $\Box$}

\noindent \textbf{Remark 3.} Another proof of Case 1 above reads as follows: Take $x\in \mathbb{Q}_p$ and $h\in p^{-N}\mathbb{Z}_p$. We want to show that $\Delta_hf(x)=0$ or, in other words, that $f(x+h)=f(x)$. Let $u\in\mathbb{Q}_p\setminus p^{-N}\mathbb{Z}_p$. Then $f(x+h+u)=f(x+h)$  because $\Delta_uf(x+h)=0$. On the other hand, $u+h \in \mathbb{Q}_p\setminus p^{-N}\mathbb{Z}_p$, so that $f(x)=f(x+u+h)$ because $\Delta_{u+h}f(x)=0$. This ends the proof. {\hfill $\Box$}

\begin{theorem}\label{dospadic} Let $f:\mathbb{Q}_p\to Y$ and let $(N_1,\cdots,N_s)\in\mathbb{Z}^s$, $(a_1,\cdots,a_s)\in\mathbb{Q}_p^s$ be such that  $\Delta_{h_1h_2\cdots h_s}f(x)=0$ for all $x\in\mathbb{Q}_p$ and all $(h_1,\cdots,h_s)\in
\left(\mathbb{Q}_p\setminus (a_1+ p^{-N_1}\mathbb{Z}_p)\right) \times \cdots\times \left(\mathbb{Q}_p\setminus (a_s+ p^{-N_s}\mathbb{Z}_p)\right)$.  Then $\Delta_{h_1h_2\cdots h_s}f(x)=0$ for all $(x,h_1,\cdots,h_s)\in\mathbb{Q}_p^{s+1}$.
\end{theorem}

\noindent \textbf{Proof.}  It is enough to apply the same arguments we used for the proof of Theorem \ref{dos}, just replacing Lemma \ref{lema} by Lemma \ref{lemapadic}. {\hfill $\Box$}


\footnotesize{J. M. Almira

Departamento de Matemáticas. Universidad de Ja\'{e}n.

E.U.P. Linares C/Alfonso X el Sabio, 28

23700 Linares (Ja\'{e}n) Spain

email: {\ttfamily jmalmira@ujaen.es}}

\end{document}